\documentclass{article}
\usepackage{amssymb,amsmath,amsthm}
\newtheorem{theorem}{Theorem}

\newtheorem{lemma}{Lemma}

\newtheorem{corollary}{Corollary}
\newtheorem{remark}{Remark}
\date{}
\numberwithin{equation}{section} \numberwithin{theorem}{section}
\numberwithin{lemma}{section} \numberwithin{corollary}{section}
\numberwithin{remark}{section} \numberwithin{proposition}{section}
\numberwithin{definition}{section}
\def\avint{\mathop{\,\rlap{--}\!\!\int}\nolimits}
\begin{document}
\newcommand{\n}{\noindent}
\newcommand{\vs}{\vskip}

\title{
Global Higher Integrability of the Gradient of the A-Laplace Equation Solution}
\author{Abdeslem Lyaghfouri\\
United Arab Emirates University\\
Department of Mathematical Sciences\\
Al Ain, Abu Dhabi, UAE}
\maketitle
\begin{abstract}
In this paper, we establish higher integrability of the gradient of the solution of the
quasilinear elliptic equation $\Delta_A u=\text{div}\left(\frac{a(|F|)}{|F|}F\right)$
in $\mathbb{R}^n$, where $\Delta_A$ is the so called $A$-Laplace operator.
\end{abstract}

\vs 0.5cm
\n Key words :   $A$-Laplace, Orlicz Sobolev space, Higher Integrability

\vs 0.5cm \noindent AMS Mathematics Subject Classification:
35J62, 35B65

\section{Introduction}\label{S:intro}

\vs 0.3cm \n We are interested in higher integrability of the solution of the following problem:

\begin{equation*} (P)\begin{cases}
&~~ u\in W^{1,A}(\mathbb{R}^n),\\
&~~\displaystyle{\Delta_A u=\text{div}\left(\Theta(F)\right)}\qquad \text{in }\mathbb{R}^n
\end{cases}\end{equation*}

\n where $n\geqslant 2$, $F=(F_1,...,F_n)\in L^{A}(\mathbb{R}^n)$, $\Delta_A u=\displaystyle{\text{div}\left(\Theta(\nabla
u)\right)}$, $\displaystyle{ \Theta(X)={{a(|X|)}\over {|X|}}X}$, and 
$\displaystyle{A(t)=\int_0^t a(s)ds}$, with $a$ a function in $C^1( (0,\infty))\cap C^0([0,\infty))$
satisfying  $a(0)=0$,  and the condition
\begin{equation}\label{1.1}
 a_0\leqslant {{ta'(t)}\over{a(t)}}\leqslant a_1    \qquad \forall t>0, \quad  a_0, a_1 \hbox{ positive
constants}.
\end{equation}

\n Without loss of generality, we shall assume that $a_0<1<a_1$.

\vs 0,5cm \n We call a solution of problem $(P)$ any
function $u\in W^{1,A}(\mathbb{R}^n)$ that satisfies
$$\int_{\mathbb{R}^n} \Theta(\nabla u).\nabla \varphi dx=\int_{\mathbb{R}^n} \Theta(F).\nabla\varphi dx\qquad
\forall\varphi\in \mathcal{D}(\mathbb{R}^n)$$

\vs 0.2cm \n We recall the definition of the Orlicz space $L^A(\mathbb{R}^n)$ and its norm (see \cite{[M]}) 
\begin{eqnarray*}
&& L^A(\mathbb{R}^n)=\left\{\,u\in L^1(\mathbb{R}^n)\,/\,\int_{\mathbb{R}^n}
A(|u(x)|)dx<\infty\,\right\}, \\
&& ||u||_A=\inf\left\{k>0\,/\,\int_{\mathbb{R}^n} A\Big({{|u(x)|}\over
k}\Big)dx\leq 1\,\right\}
\end{eqnarray*}

\n The dual space of $L^A(\mathbb{R}^n)$ is the Orlicz space $L^{\widetilde{A}}(\mathbb{R}^n)$,
where $\displaystyle{\widetilde{A}(t)=\int_{0}^{t} a^{-1}(s) ds}$ and
$a^{-1}$ is the inverse function of $a$.

\n The Orlicz-Sobolev space $W^{1,A}(\mathbb{R}^n)$ and its norm are given by
\begin{equation*}
W^{1,A}(\mathbb{R}^n)=\left\{\,u\in L^A(\mathbb{R}^n)\,/\,|\nabla u|\in
L^A(\mathbb{R}^n)\,\right\},\quad  ||u||_{1,A}=||u||_A+||\nabla u||_A.
\end{equation*}

\n Both $L^A(\mathbb{R}^n)$ and $ W^{1,A}(\mathbb{R}^n)$ are Banach, reflexive spaces.

\vs 0.5cm\n The following useful inequalities can
be easily deduced from $(1.1)$ (see  \cite{[L]} )

\begin{equation}\label{e1.2}
  {t a(t) \over 1+ a_1} \leq   A(t) \leq
  t a(t)\qquad \forall t \geq 0,
\end{equation}

\begin{equation}\label{e1.3}
 s  a(t) \leq sa(s)  +
  t a(t)\qquad \forall s,t \geq 0,
\end{equation}

\begin{equation}\label{e1.4}
  \min ( s^{a_0},s^{a_1})
   a(t) \leq a(st) \leq \max ( s^{a_0},s^{a_1})
   a(t)\qquad \forall s,t \geq 0,
\end{equation}

\begin{equation}\label{e1.5}
 \min( s^{1+a_0},s^{1+ a_1})
  {A(t) \over 1+ a_1} \leq   A(st) \leq
   {(1+ a_1)}\max ( s^{1+a_0},s^{1+ a_1})
  A(t) \qquad \forall s,t \geq 0.
\end{equation}

\n Using (1.5) and the convexity of $A$, we obtain 
\begin{eqnarray} \label{1.6}
A(s+t)&=&A\left(2.\left(\frac{s+t}{2}\right)\right)\leq(1+ a_1)2^{1+a_1} A\left(\frac{s+t}{2}\right)\nonumber \\
&\leq&(1+ a_1)2^{a_1}(A(s)+A(t))\quad \forall s,t \geq 0
\end{eqnarray}

\n We also recall the following monotonicity inequality (see \cite{[CL1]})

\begin{eqnarray}\label{e1.7}
&&\Big(\Theta(X)-\Theta(Y)\Big). (X-Y) \geqslant C(A,n)|X-Y|^2{ a\big((|X|^2+|Y|^2)^{1/2}\big)\over {(|X|^2+|Y|^2)^{1/2}}} \nonumber\\
&&\qquad \forall (X,Y)\in \mathbb{R}^{2n}\setminus\{ 0\}
\end{eqnarray}

\n where $C(A,n)$ is a positive constant depending only on $A$ and $n$.

\vs0.2cm\n There is a wide range of functions $a(t)$ satisfying (1.1). 
In particular, we observe that we have $a_0=a_1=p-1$ if and only if 
$a(t)=t^{a_0}$. In this case, we have $A(t)={{t^p}\over p}$ and
$\Delta_A$ is the $p-$Laplace operator $\Delta_p$.
We refer to \cite{[CLR]} for more examples of these functions.  

\vs 0.3cm \n Here is the main result of this paper.

\begin{theorem}\label{t1.1} Let $u\in W^{1,A}(\mathbb{R}^n)$ be a solution of $(P)$ and
assume that $F\in L^B(\mathbb{R}^n)$ for some function $B: [0,\infty)\rightarrow [0,\infty)$ 
such that $B\circ A^{-1}$ satisfies (1.1). Then we have $|\nabla u|\in L^B(\mathbb{R}^n)$
with
\[\int_{\mathbb{R}^n} B(|\nabla u|)dx\leq C\int_{\mathbb{R}^n} B(|F|) dx\]
where $C$ is a positive constant depending only on $n$, $A$, and $B$.
\end{theorem}

\vs 0,3cm\n We would like to mention that this regularity is well known for the 
Laplace equation (see \cite{[I]}). 
For the p-Laplace equation, we refer to \cite{[I]} and to \cite{[DM]} for systems. 

\vs 0,3cm
\begin{corollary}\label{c1.2}
If $u\in W^{1,A}(\mathbb{R}^n)$ is an $A-$Harmonic function in $\mathbb{R}^n$
i.e. a solution of problem $(P)$ with $F=0$, then $u$ is identically zero in $\mathbb{R}^n$.
\end{corollary}

\vs 0,3cm
\begin{remark}\label{r1.1}
i) We observe that $B=(B\circ A^{-1})\circ A$ satisfies (1.1)
as a compose of two functions satisfying the same property namely $B\circ A^{-1}$ and $A$.

\n ii) Since $B\circ A^{-1}$ satisfies (1.1), 
by using (1.5) with $t=1$ and $\displaystyle{s=\frac{1}{t}}$, we see that
there exist two positive constants $\mu>1$ and $K$ such that
\begin{equation*}
 0\leq B\circ A^{-1}\left(\frac{1}{t}\right) \leq \frac{K}{t^\mu}
    \qquad \forall t \geq 1.
\end{equation*}
Therefore the improper integral $\displaystyle{\int_1^\infty B\circ A^{-1}\left(\frac{1}{t}\right)dt}$
is convergent. This property will be used in the proof of Theorem 1.1 in Section 3.
\end{remark}

\vs 0,5cm In the sequel, we will denote by $u$ a solution of problem $(P)$.
In Section 2, we recall some well known results about
$A-$harmonic functions and establish a few Lemmas to pave the way for the proof 
of Theorem 1.1 which will be given in Section 3.

\section{Some Auxiliary Lemmas  }\label{2}

\vs 0.5cm \n Let $x_0\in \mathbb{R}^n$ and $R>0$. For each open ball $B_R(x_0)$ in $\mathbb{R}^n$ 
of center $x_0$ and radius $R$, let $v$ be the unique solution of the following problem:
\[ (P_0)\begin{cases}
& \quad v\in W^{1,A}(B_R(x_0)),\\
&\quad \Delta_A v=0\qquad \text{in }B_R(x_0), \\
& \quad v=u\qquad \text{in }\partial B_R(x_0)
\end{cases}\]

\vs 0.3cm \n First, we recall some properties of the solution of problem $(P_0)$.

\begin{lemma}\label{l2.1}
\begin{equation*}
\int_{B_R(x_0)} a(|\nabla v|) |\nabla v|  \leq
2^{a_1+2}\int_{B_R(x_0)} a(|\nabla u|) |\nabla u|dx
\end{equation*}
\end{lemma}

\n\emph{Proof}. See \cite{[CL1]}, Proof of Lemma 3.1.

\begin{remark}\label{r2.1}
Using (1.2) and Lemma 2.1, we obtain
\begin{eqnarray*}
\int_{B_R(x_0)} A(|\nabla v|)dx &\leq& \int_{B_R(x_0)} a(|\nabla v|) |\nabla v|dx\nonumber\\
&\leq& 2^{a_1+2}\int_{B_R(x_0)} a(|\nabla v|) |\nabla v|dx\nonumber\\
&\leq& (1+a_1)2^{a_1+2}\int_{B_R(x_0)} A(|\nabla u|)dx
\end{eqnarray*}
\end{remark}

\begin{lemma}\label{l2.2} There exists a positive constant $C_1=C_1(n,a_1)$ such that
\begin{equation*}
\sup_{ B_{R/2}(x_0)} A(|\nabla v|)\leq {C_1\over R^n} \int_{B_{R}(x_0)}
 A(|\nabla v|)dx
 \end{equation*}
\end{lemma}

\n\emph{Proof}. See \cite{[L]}, Lemma 5.1.

\vs 0.2cm\n For each function $f: \mathbb{R}^n \rightarrow \mathbb{R}~(\mathbb{R}^n)$, let
$\displaystyle{(f)_{x_0,r}=\frac{1}{|B_r(x_0)|}\int_{B_{r}(x_0)} f dx=\avint_{B_{r}(x_0)} f dx }$
and $\displaystyle{(f)_{r}=(f)_{0,r} }$. Then we have the following property of the function $v$.

\begin{lemma}\label{l2.3} There exist two positive constants $\alpha=\alpha(n,a_1)<1$ and $C_2=C_2(n,a_1)$ such that
we have for any $r\in(0,R)$
\begin{equation*}
\avint_{B_{r}(x_0)} A(|\nabla v-(\nabla v)_{x_0,r}|)dx\leq C_2\left({r\over R}\right)^\alpha\avint_{B_{R}(x_0)}
 A(|\nabla v-(\nabla v)_{x_0,R}|)dx
 \end{equation*}
\end{lemma}

\n\emph{Proof}. See \cite{[L]}, Lemma 5.1.

\begin{remark}\label{r2.2}
Using (1.6), Remark 2.1, and Jensen's inequality, we obtain
\begin{eqnarray}\label{e2.1}
\avint_{B_{R}(x_0)}A(|\nabla v-(\nabla v)_{x_0,R}|)dx &\leq& \avint_{B_{R}(x_0)}
 A(|\nabla v|+|(\nabla v)_{x_0,R}|)dx\nonumber\\
&\leq& (1+a_1)2^{a_1}\left(\avint_{B_{R}(x_0)}
A(|\nabla v|)dx+A(|(\nabla v)_{x_0,R}|) \right)\nonumber\\
&=& (1+a_1)2^{a_1}\left(\avint_{B_{R}(x_0)}
A(|\nabla v|)dx+A\left(\avint_{B_{R}(x_0)}|\nabla v|dx\right) \right)\nonumber\\
&\leq& (1+a_1)2^{a_1}\left(\avint_{B_{R}(x_0)}
A(|\nabla v|)dx+\avint_{B_{R}(x_0)}A(|\nabla v|)dx \right)\nonumber\\
&=& (1+a_1)2^{a_1+1}\avint_{B_{R}(x_0)}A(|\nabla v|)dx\nonumber\\
&\leq&(1+a_1)^22^{2a_1+3}\avint_{B_{R}(x_0)}A(|\nabla u|)dx
\end{eqnarray}

\n Combining (2.1) and Lemma 2.3, we obtain
\begin{equation*}
\avint_{B_{r}(x_0)} A(|\nabla v-(\nabla v)_{x_0,r}|)dx\leq
C_2(1+a_1)^22^{2a_1+3}\left({r\over R}\right)^\alpha\avint_{B_{R}(x_0)}A(|\nabla u|)dx
\end{equation*}
\end{remark}

\vs 0.3cm \n The next lemma is the key tool in the proof of Theorem 1.1.

\begin{lemma}\label{l2.4}
There exist two positive constants $\gamma=\gamma(n,a_0,a_1)$ and $m=m(\alpha,n,a_0,a_1)$
such that for each $\delta\in(0,1)$, we have for any $x_0\in \mathbb{R}^n$, $R>0$ and $r\in(0,R)$:

 \begin{eqnarray*}
&&\avint_{B_{r}(x_0)} |A(|\nabla u|)-(A(|\nabla u|))_r|dx
~\leq~\frac{\gamma}{\delta^{a_1(1+a_1)}}\left({R\over r}\right)^m\avint_{B_{R}(x_0)}A(|F|)dx\nonumber\\
&&\quad+\gamma\left(\delta^{a_0+1}+\frac{1}{\delta^{a_1(1+a_1)}}\left({r\over R}\right)^\alpha\right)\avint_{B_{R}(x_0)}A(|\nabla u|)dx
\end{eqnarray*}
\end{lemma}

\vs 0.3cm\n The proof of Lemma 2.4 requires several lemmas.

\begin{lemma}\label{l2.5}  Let $G\,:\, \mathbb{R}^{2n}\setminus\{0\}\,\rightarrow\,
\mathbb{R}$ defined by

$$G(\xi,\zeta)=\int_0^1 {{a(|\theta_t|)}\over   {|\theta_t|}}dt,\quad \theta_t=t\xi+(1-t)\zeta.$$

\n Then there exists two positive constants $c_{a_1,n}$ depending only on $n$ and $a_1$, and $C_{a_0,n}$
depending only on $n$ and $a_0$ such that:
\[\forall (\xi,\zeta)\in\mathbb{R}^{2n}\setminus\{0\},\quad c_{a_1,n} {{a(|\xi|+|\zeta|)}\over   {|\xi|+|\zeta|}}\leq
G(\xi,\zeta) \leq C_{a_0,n}
{{a(|\xi|+|\zeta|)}\over{|\xi|+|\zeta|}}\]
\end{lemma}

\vs 0,3cm \n The proof of Lemma 2.5 is based on the following
lemma proved in \cite{[AS]} for $n=2$ and whose proof extends easily to $n\geq3$.

\begin{lemma}\label{l2.6}  Let $n\geqslant 2$, $p>1$ and let $F_p\,:\, \mathbb{R}^{2n}\setminus\{0\}\,\rightarrow\,
\mathbb{R}$ defined by
\[F_p(\xi,\zeta)=\int_0^1 |t\xi+(1-t)\zeta|^{p-2}dt\]
Then there exists two positive constants $c(p,n)<C(p,n)$ depending only on $p$ and $n$  such that:
\[\forall (\xi,\zeta)\in\mathbb{R}^{2n}\setminus\{0\},\quad c(p,n)\big(|\xi|^2+|\zeta|^2\big)^{{p-2}\over 2}\leq
F_p(\xi,\zeta) \leq C(p,n)\big(|\xi|^2+|\zeta|^2\big)^{{p-2}\over 2}\]
\end{lemma}

\vs 0,3cm \n \emph{Proof of Lemma 2.5.} For
$(\xi,\zeta)\in\mathbb{R}^{2n}\setminus\{0\}$, we set

$$X_0={\xi\over{(|\xi|^2+|\zeta|^2)^{1/2}}},\quad
X_1={\zeta\over{(|\xi|^2+|\zeta|^2)^{1/2}}}.$$

\n It follows that
\begin{equation*}
  \theta_t=(|\xi|^2+|\zeta|^2)^{1/2}(tX_0+(1-t)X_1),\quad
|\theta_t|=(|\xi|^2+|\zeta|^2)^{1/2}|tX_0+(1-t)X_1|.
\end{equation*}

\n Then we have by inequality (1.4)

\begin{eqnarray*}
 & a(|\theta_t|)\leq
\max\big(|tX_0+(1-t)X_1|^{a_0},|tX_0+(1-t)X_1|^{a_1}\big)a\big((|\xi|^2+|\zeta|^2)^{1/2}\big)   \\
&\min\big(|tX_0+(1-t)X_1|^{a_0},|tX_0+(1-t)X_1|^{a_1}\big)a\big((|\xi|^2+|\zeta|^2)^{1/2}\big)
\leq a(|\theta_t|)
\end{eqnarray*}

\n Note that since $|tX_0+(1-t)X_1|\leq t|X_0|+(1-t)|X_1|\leq t+(1-t)=1$
and $ 0< a_0 \leq a_1$, we have $|tX_0+(1-t)X_1|^{a_1} \leq |tX_0+(1-t)X_1|^{a_0}.$
Hence we get
\begin{equation*}
|tX_0+(1-t)X_1|^{a_1}a\big((|\xi|^2+|\zeta|^2)^{1/2}\big)
\leq  a(|\theta_t|)\leq
|tX_0+(1-t)X_1|^{a_0}a\big((|\xi|^2+|\zeta|^2)^{1/2}\big) 
\end{equation*}

\n which leads to

$$F_{a_1+1}(X_0,X_1){{a\big((|\xi|^2+|\zeta|^2)^{1/2}\big)}\over
{(|\xi|^2+|\zeta|^2)^{1/2}}}\leq G(\xi,\zeta) \leq
F_{a_0+1}(X_0,X_1){{a\big((|\xi|^2+|\zeta|^2)^{1/2}\big)}\over
{(|\xi|^2+|\zeta|^2)^{1/2}}}.$$

\n Now, if we apply Lemma 2.6 with $p=a_0+1$ and $p=a_1+1$, we get
since $|X_0|^2+|X_1|^2=1$

\begin{eqnarray*}
 &F_{a_0+1}(X_0,X_1)\leq
C(a_0+1,n)\big(|X_0|^2+|X_1|^2\big)^{{a_0-1}\over 2}=C(a_0+1,n)
 \\
&F_{a_1+1}(X_0,X_1)\geq
c_{a_1+1,n}\big(|X_0|^2+|X_1|^2\big)^{{a_1-1}\over 2}=c(a_1+1,n)
\end{eqnarray*}

\n Finally, we obtain

$$c(a_1+1,n){{a\big((|\xi|^2+|\zeta|^2)^{1/2}\big)}\over
{(|\xi|^2+|\zeta|^2)^{1/2}}}\leqslant G(\xi,\zeta) \leqslant
C(a_0+1,n){{a\big((|\xi|^2+|\zeta|^2)^{1/2}\big)}\over
{(|\xi|^2+|\zeta|^2)^{1/2}}}$$

\n Since the two norms $|\xi|+|\zeta|$ and $(|\xi|^2+|\zeta|^2)^{1/2}$
are equivalent, the lemma follows by using (1.4)

\qed

\begin{lemma}\label{l2.7}
For any $X, Y\in \mathbb{R}^n$, we have:
\begin{equation*}
A(|X|)\geq A(|Y|)+\left<\Theta(Y),X-Y\right>
\end{equation*}
\end{lemma}

\n\emph{Proof}. Let $X, Y\in \mathbb{R}^n$. First, we have
\begin{eqnarray*}
A(|X|)-A(|Y|)&=&\int_0^1 \frac{d}{dt}[A(|tX+(1-t)Y|)]dt\nonumber\\
&=& \int_0^1 A'(|tX+(1-t)Y|).\frac{<tX+(1-t)Y,X-Y>}{|tX+(1-t)Y|}dt \nonumber\\
&=& \int_0^1 a(|tX+(1-t)Y|).\frac{<tX+(1-t)Y,X-Y>}{|tX+(1-t)Y|}dt \nonumber\\
&=& \int_0^1 \left<\Theta(tX+(1-t)Y),X-Y\right>dt \nonumber\\
\end{eqnarray*}

\n Next, by using (1.7), this leads to
\begin{eqnarray*}
A(|X|)-A(|Y|)&=& \int_0^1 \frac{1}{t}\left<\Theta(tX+(1-t)Y),t(X-Y)\right>dt\\
&\geq& \int_0^1 \frac{1}{t}\left<\Theta(Y),t(X-Y)\right>dt\nonumber\\
&=& \left<\Theta(Y),X-Y\right>
\end{eqnarray*}

\n which achieves the proof.
\qed

\begin{lemma}\label{l2.8}
For any $\delta\in(0,1)$, $X, Y\in \mathbb{R}^n$, we have:
\begin{equation*}
|A(|X|)-A(|Y|)|\leq\frac{C_3}{\delta^{a_1(1+a_1)}}A(|X-Y|)+C_4\delta^{a_0+1}(A(|X|)+A(|Y|))
 \end{equation*}
where $C_3=(1+a_1)C_{a_0,n}$ and $C_4=(1+a_1)2^{a_1}C_3$
\end{lemma}

\n\emph{Proof}. Let $\delta\in(0,1)$, $X, Y\in \mathbb{R}^n$. First, we have
\begin{eqnarray*}
|A(|X|)-A(|Y|)|&=& \left|\int_0^1 \frac{d}{dt}[A(|tX+(1-t)Y|)]dt\right|\nonumber\\
&=& \left|\int_0^1 A'(|tX+(1-t)Y|).\frac{<tX+(1-t)Y,X-Y>}{|tX+(1-t)Y|}dt\right| \nonumber\\
&=& \left|\int_0^1 a(|tX+(1-t)Y|).\frac{<tX+(1-t)Y,X-Y>}{|tX+(1-t)Y|}dt\right| \nonumber\\
&\leq&|X-Y|.(|X|+|Y|).\int_0^1 \frac{a(|tX+(1-t)Y|)}{|tX+(1-t)Y|}dt
\end{eqnarray*}

\n Next, we get by using Lemma 2.5
\begin{eqnarray}\label{e2.2}
|A(|X|)-A(|Y|)|&\leq& C_{a_0,n}|X-Y|.(|X|+|Y|). {{a(|X|+|Y|)}\over
{|X|+|Y|}}\nonumber\\
&=& C_{a_0,n}|X-Y|.a(|X|+|Y|)
\end{eqnarray}

\n We observe that we have by (1.2)-(1.4), for $s, t\geq 0$ since $\delta\in(0,1)$
\begin{eqnarray}\label{e2.3}
s a(t)&=&\frac{s}{\delta^{a_1}}.\delta^{a_1}a(t)\leq \frac{s}{\delta^{a_1}}.a(\delta t)
\leq \frac{s}{\delta^{a_1}}.a\left(\frac{s}{\delta^{a_1}}.\right)  +\delta t a(\delta t)\nonumber\\
&\leq& \frac{s}{\delta^{a_1}}.\frac{1}{(\delta^{a_1})^{a_1}}a(s)  +\delta t \delta^{a_0}a(t)=
\frac{sa(s)}{\delta^{a_1(1+a_1)}}+\delta^{a_0+1}ta(t)\nonumber\\
&\leq& \frac{1+a_1}{\delta^{a_1(1+a_1)}}A(s)+(1+a_1)\delta^{a_0+1}A(t)
\end{eqnarray}

\n Using (2.2) and (2.3), with $s=|X-Y|$ and $t=|X|+|Y|$, we get 
\begin{eqnarray}\label{e2.4}
&&|A(|X|)-A(|Y|)|\leq\frac{(1+a_1)C_{a_0,n}}{\delta^{a_1(1+a_1)}}A(|X-Y|)+(1+a_1)C_{a_0,n}\delta^{a_0+1}A(|X|+|Y|)\nonumber\\
\end{eqnarray}

\n Using (1.6), we get from (2.4)
\begin{eqnarray*}
|A(|X|)-A(|Y|)|\leq\frac{C_3}{\delta^{a_1(1+a_1)}}A(|X-Y|)+(1+a_1)^22^{a_1}C_{a_0,n}\delta^{a_0+1}(A(|X|)+A(|Y|))
\end{eqnarray*}

\n Hence, the lemma follows.
\qed

\begin{lemma}\label{l2.9}
There exist two positive constants $\displaystyle{C_5=\frac{(1+a_1)2^{\frac{a_1}{2}}}{C(A,n)} }$
and 

$\displaystyle{C_6=(1+a_1)2^{a_1}\left(1+(1+a_1)2^{a_1+2}\right)\left(C_5+\frac{1+a_1}{4}\right) }$
such that we have for each $\eta\in(0,1)$ , $x_0\in \mathbb{R}^n$, $R>0$ and $r\in(0,R)$:
\begin{eqnarray*}
&&\avint_{B_{r}(x_0)} A(|\nabla u-\nabla v|)dx\leq \left({R\over r}\right)^n\Bigl(\frac{C_5}{\eta^{1+a_1(1+a_1)}}\avint_{B_{R}(x_0)}
 A(|F|)dx\nonumber\\
&&\quad+C_6\eta^{a_0}\avint_{B_{R}(x_0)}A(|\nabla u|)dx\Bigl)
 \end{eqnarray*}
\end{lemma}

\vs 0.3cm\n\emph{Proof}. We observe that by translation, it is enough to prove the lemma for $x_0=0$.
Using $w=(u-v)\chi_{B_R}$ as a test function for problems $(P)$ and $(P_0)$ and subtracting
the two equations and using (2.3), we get for $\eta\in(0,1)$:
\begin{eqnarray}\label{e2.5}
&&\int_{B_R}\big(\Theta(\nabla u)-\Theta(\nabla v)\big)(\nabla u-\nabla v)dx=\int_{B_R}\frac{a(|F|)}{|F|}F.\nabla w dx 
\leq\int_{B_R} a(|F|).|\nabla w| dx\nonumber\\
&&\qquad\leq~ \frac{1+a_1}{\eta^{a_1(1+a_1)}}\int_{B_R}A(|F|)dx+(1+a_1)\eta^{a_0+1}\int_{B_R}A(|\nabla w|)dx
\end{eqnarray}

\n Using (2.5) and (1.7), we get
\begin{eqnarray*}
\frac{C(A,n)}{2^{\frac{a_1}{2}}}\int_{B_R}{{a(|\nabla u|+|\nabla v|)}\over{|\nabla u|+|\nabla v|}}.|\nabla w|^2dx&\leq& \frac{1+a_1}{\eta^{a_1(1+a_1)}}\int_{B_R}A(|F|)dx\\
&&+(1+a_1)\eta^{a_0+1}\int_{B_R}A(|\nabla w|)dx
\end{eqnarray*}

\n or
\begin{eqnarray}\label{e2.6}
&&\int_{B_R}{{a(|\nabla u|+|\nabla v|)}\over{|\nabla u|+|\nabla v|}}.|\nabla w|^2dx~\leq~ \frac{C_5}{\eta^{a_1(1+a_1)}}\int_{B_R}A(|F|)dx+C_5\eta^{a_0+1}\int_{B_R}A(|\nabla w|)dx\nonumber\\
&&
\end{eqnarray}

\n By (1.2), we can write
\begin{equation}\label{e2.7}
\int_{B_r}A(|\nabla w|)dx\leq\int_{B_r}|\nabla w|a(|\nabla w|)dx=\int_{E_{1,r}}|\nabla w|a(|\nabla w|)dx+\int_{E_{2,r}}|\nabla w|a(|\nabla w|)dx
\end{equation}

\n where $E_{1,r}=B_r\cap\{~(|\nabla u|+|\nabla v|)a(|\nabla w|)\leq |\nabla w|a(|\nabla u|+|\nabla v|)~\}$
and $E_{2,r}=B_r\cap\{~(|\nabla u|+|\nabla v|)a(|\nabla w|)>|\nabla w|a(|\nabla u|+|\nabla v|)~\}$.

\n From the definition of $E_{1,r}$, we see that
\begin{eqnarray}\label{e2.8}
\int_{E_{1,r}}|\nabla w|a(|\nabla w|)dx\leq\int_{E_{1,r}}{{a(|\nabla u|+|\nabla v|)}\over{|\nabla u|+|\nabla v|}}.|\nabla w|^2dx
\end{eqnarray}
\n Using the monotonicity of $a$, the definition of $E_{2,r}$, and Young's inequality, we get
\begin{eqnarray}\label{e2.9}
&&\int_{E_{2,r}}|\nabla w|a(|\nabla w|)dx~\leq~\int_{E_{2,r}}a(|\nabla u|+|\nabla v|).|\nabla w|dx\nonumber\\
&&\quad\leq\int_{E_{2,r}}\left(\frac{a(|\nabla u|+|\nabla v|)}{|\nabla u|+|\nabla v|}.|\nabla w|^2\right)^{\frac{1}{2}}
.\left((|\nabla u|+|\nabla v|)a(|\nabla u|+|\nabla v|)\right)^{\frac{1}{2}}dx\nonumber\\
&&\quad\leq\frac{1}{\eta}\int_{E_{2,r}}{{a(|\nabla u|+|\nabla v|)}\over{|\nabla u|+|\nabla v|}}.|\nabla w|^2dx
+\frac{\eta}{4}\int_{E_{2,r}}(|\nabla u|+|\nabla v|)a(|\nabla u|+|\nabla v|)dx\nonumber\\
&&\quad\leq\frac{1}{\eta}\int_{E_{2,r}}{{a(|\nabla u|+|\nabla v|)}\over{|\nabla u|+|\nabla v|}}.|\nabla w|^2dx
+\frac{(1+a_1)\eta}{4}\int_{B_r}A(|\nabla u|+|\nabla v|)dx\nonumber\\
\end{eqnarray}

\n Combing (2.6), (2.7), (2.8), and (2.9), and using (1.6) and Remark 2.1, we arrive at
\begin{eqnarray*}
&&\int_{B_r}A(|\nabla w|)dx~\leq~\frac{1}{\eta}\int_{B_r}{{a(|\nabla u|+|\nabla v|)}\over{|\nabla u|+|\nabla v|}}.|\nabla w|^2dx
+\frac{(1+a_1)\eta}{4}\int_{B_R}A(|\nabla u|+|\nabla v|)dx\nonumber\\
&&\quad\leq \frac{C_5}{\eta^{a_1(1+a_1)+1}}\int_{B_R}A(|F|)dx+C_5\eta^{a_0}\int_{B_R}A(|\nabla w|)dx
+\frac{(1+a_1)\eta}{4}\int_{B_R}A(|\nabla u|+|\nabla v|)dx\nonumber\\
&&\quad\leq \frac{C_5}{\eta^{a_1(1+a_1)+1}}\int_{B_R}A(|F|)dx+\left(C_5\eta^{a_0}
+\frac{(1+a_1)\eta}{4}\right)\int_{B_R}A(|\nabla u|+|\nabla v|)dx\nonumber\\
&&\quad\leq (1+a_1)2^{a_1}\left(C_5\eta^{a_0}
+\frac{(1+a_1)\eta}{4}\right)\left(\int_{B_R}A(|\nabla u|)dx+\int_{B_R}A(|\nabla v|)dx\right)\nonumber\\
&&\quad+\frac{C_5}{\eta^{a_1(1+a_1)+1}}\int_{B_R}A(|F|)dx\nonumber\\
&&\quad\leq (1+a_1)2^{a_1}\left(C_5
+\frac{1+a_1}{4}\right)(1+(1+a_1)2^{a_1+2})\eta^{a_0}\int_{B_R}A(|\nabla u|)dx\nonumber\\
&&\quad+\frac{C_5}{\eta^{a_1(1+a_1)+1}}\int_{B_R}A(|F|)dx
\end{eqnarray*}

\n Hence the lemma follows.
\qed

\vs 0.3cm\n\emph{Proof of Lemma 2.4}. Let $\delta\in(0,1)$, $R>0$ and $r\in(0,R)$. First, we have by
Lemma 2.8 used for $X=\nabla u(x)$ and $Y=\nabla u(y)$

\begin{eqnarray*}
&&\avint_{B_{r}} |A(|\nabla u|)-(A(|\nabla u|))_r|)dx~=~\avint_{B_{r}}\left|\avint_{B_{r}}(A(|\nabla u(x)|)-A(|\nabla u(y)|))dy\right|dx\nonumber\\
&&\quad\leq\avint_{B_{r}\times B_{r}} |A(|\nabla u(x)|)-A(|\nabla u(y)|)|dxdy\nonumber\\
&&\quad\leq\frac{C_3}{\delta^{a_1(1+a_1)}}\avint_{B_{r}\times B_{r}} A(|\nabla u(x)-\nabla u(y)|) dxdy\nonumber\\
&&\qquad+C_4\delta^{a_0+1}\left(\avint_{B_{r}\times B_{r}}A(|\nabla u(x)|)dxdy+\avint_{B_{r}\times B_{r}}A(|\nabla u(y)|)dxdy\right)\nonumber\\
&&\quad=\frac{C_3}{\delta^{a_1(1+a_1)}}\avint_{B_{r}\times B_{r}} A(|\nabla u(x)-\nabla u(y)|) dxdy\nonumber\\
&&\qquad+2C_4\delta^{a_0+1}\avint_{B_{r}}A(|\nabla u(x)|)dx\nonumber\\
&&\quad\leq\frac{C_3}{\delta^{a_1(1+a_1)}}\avint_{B_{r}\times B_{r}}A(|\nabla u(x)-(\nabla u)_r|+|\nabla u(y)-(\nabla u)_r|)dxdy\nonumber\\
&&\qquad+2C_4\delta^{a_0+1}\avint_{B_{r}}A(|\nabla u|)dx\nonumber
\end{eqnarray*}

\n This leads by (1.6) to
\begin{eqnarray}\label{e2.10}
&&\avint_{B_{r}} |A(|\nabla u|)-(A(|\nabla u|))_r|dx\nonumber\\
&&\quad\leq\frac{C_3(1+ a_1)2^{a_1}}{\delta^{a_1(1+a_1)}}\left(\avint_{B_{r}\times B_{r}} A(|\nabla u(x)-(\nabla u)_r|)dxdy
+\avint_{B_{r}\times B_{r}} A(|\nabla u(y)-(\nabla u)_r|)dxdy\right)\nonumber\\
&&\qquad+2C_4\delta^{a_0+1}\avint_{B_{r}}A(|\nabla u|)dx\nonumber\\
&&\quad=\frac{2C_3(1+ a_1)2^{a_1}}{\delta^{a_1(1+a_1)}}\avint_{B_{r}} A(|\nabla u-(\nabla u)_r|)dx
+2C_4\delta^{a_0+1}\avint_{B_{r}}A(|\nabla u|)dx
\end{eqnarray}

\n Using (1.6) again, we get
\begin{eqnarray}\label{e2.11}
&&\avint_{B_{r}} A(|\nabla u-(\nabla u)_r|)dx~\leq~\avint_{B_{r}} A(|\nabla u-\nabla v|+|\nabla v-(\nabla v)_r|
+|(\nabla u)_r-(\nabla v)_r|)dx\nonumber\\
&&\quad\leq(1+ a_1)2^{a_1}\left(\avint_{B_{r}} A(|\nabla u-\nabla v|)dx
+\avint_{B_{r}} A(|\nabla v-(\nabla v)_r|+|(\nabla u)_r-(\nabla v)_r|)dx\right)\nonumber\\
&&\quad\leq(1+ a_1)2^{a_1}\avint_{B_{r}} A(|\nabla u-\nabla v|)dx\nonumber\\
&&\quad
+(1+ a_1)^22^{2a_1}\left(\avint_{B_{r}} A(|\nabla v-(\nabla v)_r|)dx+\avint_{B_{r}} A(|(\nabla u)_r-(\nabla v)_r|)dx\right)\nonumber\\
\end{eqnarray}

\n Using Jensen's inequality, we derive
\begin{eqnarray}\label{e2.12}
&&\avint_{B_{r}} A(|(\nabla u)_r-(\nabla v)_r|)dx=A\left(\left|\avint_{B_{r}}(\nabla u-\nabla v)dx\right|\right)
\leq\avint_{B_{r}} A(|\nabla u-\nabla v|)dx\nonumber\\
&&
\end{eqnarray}

\n From (2.11) and (2.12), we obtain
\begin{eqnarray}\label{e2.13}
&&\avint_{B_{r}} A(|\nabla u-(\nabla u)_r|)dx~\leq~(1+ a_1)2^{a_1}(1+(1+ a_1)2^{a_1})\avint_{B_{r}} A(|\nabla u-\nabla v|)dx\nonumber\\
&&\quad +(1+ a_1)^22^{2a_1}\avint_{B_{r}} A(|\nabla v-(\nabla v)_r|)dx
\end{eqnarray}

\n Combining (2.10) and (2.13), we obtain
\begin{eqnarray}\label{e2.14}
&&\avint_{B_{r}} |A(|\nabla u|)-(A(|\nabla u|))_r|dx\leq~\frac{C_7}{\delta^{a_1(1+a_1)}}\avint_{B_{r}} A(|\nabla u-\nabla v|)dx\nonumber\\
&&\quad+\frac{C_8}{\delta^{a_1(1+a_1)}}\avint_{B_{r}} A(|\nabla v-(\nabla v)_r|)dx+2C_4\delta^{a_0+1}\avint_{B_{r}}A(|\nabla u|)dx\nonumber\\
\end{eqnarray}

\n where $\displaystyle{C_7=2C_3(1+ a_1)^22^{2a_1}(1+(1+ a_1)2^{a_1}) }$
and $\displaystyle{C_8=2C_3(1+ a_1)^32^{3a_1} }$.

\n Finally, by taking into account Remark 2.2 and Lemma 2.8, we obtain from (2.14) for any
$\eta\in(0,1)$ 
\begin{eqnarray*}
&&\avint_{B_{r}} |A(|\nabla u|)-(A(|\nabla u|))_r|dx\nonumber\\
&&\quad\leq~\frac{C_7}{\delta^{a_1(1+a_1)}}\left({R\over r}\right)^n\Bigl(\frac{C_5}{\eta^{1+a_1(1+a_1)}}\avint_{B_{R}}
 A(|F|)dx+C_6\eta^{a_0}\avint_{B_{R}(x_0)}A(|\nabla u|)dx\Bigl)\nonumber\\
&&\quad+\frac{C_8C_2}{\delta^{a_1(1+a_1)}}(1+a_1)^22^{2a_1+3}\left({r\over R}\right)^\alpha\avint_{B_{R}}A(|\nabla u|)dx
+2C_4\delta^{a_0+1}\avint_{B_{r}}A(|\nabla u|)dx
\end{eqnarray*}

\n To conclude the proof, we let $\gamma=\max(C_5C_7,C_6C_7,C_2C_8(1+a_1)^22^{2a_1+3},2C_4)$,
$\displaystyle{m=n+\frac{(\alpha+n)(1+a_1(1+a_1))}{a_0}  }$, 
and choose $\eta$ such that $\displaystyle{\eta^{a_0}=\left({r\over R}\right)^{\alpha+n} }$
or $\displaystyle{\eta=\left({r\over R}\right)^{\frac{\alpha+n}{a_0}}  }$.
Then we obtain
\begin{eqnarray*}
&&\avint_{B_{r}} |A(|\nabla u|)-(A(|\nabla u|))_r|dx
~\leq~\frac{\gamma}{\delta^{a_1(1+a_1)}}\left({R\over r}\right)^m\avint_{B_{R}}A(|F|)dx\nonumber\\
&&\quad+\gamma\left(\delta^{a_0+1}+\frac{1}{\delta^{a_1(1+a_1)}}\left({r\over R}\right)^\alpha\right)\avint_{B_{R}}A(|\nabla u|)dx
\end{eqnarray*}
\qed

\section{Proof of Theorem 1.1}\label{3}

\vs 0.3cm\n This section is devoted to the proof of Theorem 1.1. First, we recall 
that for each function $f\in L^1(\mathbb{R}^n)$, 
the Hardy-Littlewood maximal function associated with $f$ is given by
\[M[f](x_0)=\sup_{r>0} \avint_{B_r(x_0)}|f(x)|dx,\quad x_0\in \mathbb{R}^n \] 

\n and the sharp maximal function associated with $f$ is defined by
\[f^{\sharp}(x_0)=\sup_{r>0} \avint_{B_r(x_0)}|f-(f)_{x_0,r}|dx,\quad x_0\in \mathbb{R}^n\] 

\vs 0.3cm\n\emph{Proof of Theorem 1.1}. We deduce from Lemma 2.4 that we have for every $\delta\in(0,1)$,
$r>0$, $\displaystyle{R=\frac{r}{\delta^{\frac{a_1(1+a_1)+a_0+1}{\alpha}} } }$,
$\displaystyle{\kappa=a_1(1+a_1)+\frac{m(a_1(1+a_1)+a_0+1)}{\alpha} }$, and $x_0\in \mathbb{R}^n$
\begin{eqnarray*}
&&\avint_{B_{r}(x_0)} |A(|\nabla u|)-(A(|\nabla u|))_r|)dx
~\leq~\frac{\gamma}{\delta^{\kappa}}\avint_{B_{R}x_0)}A(|F|)dx\nonumber\\
&&\qquad+2\gamma \delta^{a_0+1}\avint_{B_{R}x_0)}A(|\nabla u|)dx
\end{eqnarray*}
\n which leads to
\begin{eqnarray}\label{e3.1}
(A(|\nabla u|))^{\sharp}(x_0)~\leq~\frac{\gamma}{\delta^{\kappa}}M[A(|F|)](x_0)
+2\gamma \delta^{a_0+1} M[A(|\nabla u|)](x_0)
\end{eqnarray}

\n Now, let $B$ be a function such that $B\circ A^{-1}$ satisfies (1.1) and
assume that $F\in C_0^\infty(\mathbb{R}^n)$ and $|\nabla u|\in L^B(\mathbb{R}^n)$. 
Then using (1.6), we obtain from (3.1) (see \cite{[F]}) for 
$\gamma'=\gamma (1+c_1)2^{c_1}$, where $c_1$ is the equivalent to the constant 
$a_1$ for the function $B\circ A^{-1}$ 
\begin{eqnarray}\label{e3.2}
\int_{\mathbb{R}^n} B(|\nabla u|)dx &\leq& \int_{\mathbb{R}^n}B\left(A^{-1}(A(|\nabla u|))^{\sharp}(x)\right)dx \nonumber\\
&\leq&2\gamma'\delta^{1+a_0} \int_{\mathbb{R}^n}B\left(A^{-1}M[A(|\nabla u|)](x)\right)dx\nonumber\\
&&\quad+\frac{\gamma'}{\delta^{\kappa}}\int_{\mathbb{R}^n}B\left(A^{-1}M[A(|F|)](x)\right)dx
\end{eqnarray}

\n By the Hardy-Littlewood maximal theorem (see \cite{[G]}), we have for some positive
constant $C_9$ depending only upon $n, A$ and $B$
\begin{eqnarray}\label{e3.3-4}
\int_{\mathbb{R}^n}B\left(A^{-1}M[A(|F|)](x)\right)dx &\leq& C_9\int_{\mathbb{R}^n}B(|F|)dx\\
\int_{\mathbb{R}^n}B\left(A^{-1}M[A(|\nabla u|)](x)\right)dx &\leq& C_9\int_{\mathbb{R}^n}B(|\nabla u|)dx
\end{eqnarray}

\n We deduce then from (3.2), (3.3) and (3.4) that
\begin{equation}\label{e3.5}
\int_{\mathbb{R}^n} B(|\nabla u|)dx 
\leq2\gamma' C_9\delta^{1+a_0} \int_{\mathbb{R}^n}B(|\nabla u|)dx
+\frac{\gamma' C_9}{\delta^{\kappa}}\int_{\mathbb{R}^n}B(|F|)dx
\end{equation}

\n Now, if we choose $\delta$ such that $\delta\in\left(0, \min\left(1,\left(2\gamma' C_9\right)^{-\frac{1}{1+a_0}}\right)  \right)$,
we get
\begin{equation}\label{e3.6}
\int_{\mathbb{R}^n} B(|\nabla u|)dx
\leq \frac{\gamma' C_9}{\delta^{\kappa}(1-2\gamma' C_9\delta^{1+a_0})}\int_{\mathbb{R}^n}B(|F|)dx
\end{equation}

\n This completes the proof when $F\in C_0^\infty(\mathbb{R}^n)$ and $|\nabla u|\in L^B(\mathbb{R}^n)$.
Next, we will establish it when $F\in L^{B}(\mathbb{R}^n)$. To do that, we consider a
sequence $(F_k)_k$ of vector functions in $C_0^\infty(\mathbb{R}^n)$ that converges to $F$ 
in $L^{B}(\mathbb{R}^n)$. We denote by $(u_k)_k$ the sequence of unique 
solutions of the following problem
\begin{equation*} (P_k)\begin{cases}
& \quad u_k\in W^{1,A}(\mathbb{R}^n),\\
&\quad\displaystyle{\text{div}\big(\Theta(\nabla
u_k)\big)=\text{div}\left(\Theta(F_k)\right)}\qquad \text{in }\mathbb{R}^n
\end{cases}\end{equation*}

\n Since $F_k$ has compact support in $\mathbb{R}^n$, there exists $l_k>0$ such that
$\text{suppt}(F_k)\subset B_{l_k}$, and therefore we have for all $k$

\begin{equation}\label{e3.7} 
\text{div}\big(\Theta(\nabla u_k)\big)=0\quad \text{in }\mathbb{R}^n\setminus B_{l_k}
\end{equation}

\n We will prove that $|\nabla u_k|\in L^B(\mathbb{R}^n)$. Since $u_k\in C^1(\mathbb{R}^n)$,
is enough to show that $\displaystyle{\int_{\{|x|>2l_k\}} B(|\nabla u_k(x)|)dx<\infty}$.
We observe that because of (3.7), we can apply Lemma 2.3 in $B_{|x|-l_k}(x)$ for each $x\in\mathbb{R}^n\setminus B_{2l_k}$
\begin{equation}\label{e3.8}
A(|\nabla u_k(x)|)\leq {C_1\over (|x|-l_k)^n} \int_{B_{|x|-l_k}(x)}
 A(|\nabla u_k|)dy
 \end{equation}
 
 \n Since $B\circ A^{-1}$ satisfies (1.1), there exists by (1.5) 
two positive constants $\mu$ and $K$ such that
\begin{equation}\label{e3.9}
 0\leq B\circ A^{-1}(st) \leq K s^\mu B\circ A^{-1}(t)
    \qquad \forall s, t \geq 0.
\end{equation}

\n Using (3.8) and (3.9), and taking into account Remark 1.1, we get
\begin{eqnarray*}
&&\int_{\{|x|>2l_k\}} B(|\nabla u_k(x)|)dx\nonumber\\
&&\quad\leq K\left(C_1\int_{B_{|x|-l_k}(x)}A(|\nabla u_k|)dy\right)^\mu. \int_{\{|x|>2l_k\}}B\circ A^{-1}\left({{1}\over (|x|-l_k)^n}\right)dx \nonumber\\
&&\quad\leq K\omega_n\left(C_1\int_{\mathbb{R}^n}A(|\nabla u_k|)dy\right)^\mu. \int_{2l_k}^\infty r^{n-1}B\circ A^{-1}\left({{1}\over (r-l_k)^n}\right)dr \nonumber\\
&&\quad\leq K\omega_n2^{n-2}\left(C_1\int_{\mathbb{R}^n}A(|\nabla u_k|)dy\right)^\mu. \int_{2l_k}^\infty (r-l_k)^{n-1}B\circ A^{-1}\left({{1}\over (r-l_k)^n}\right)dr \nonumber\\
&&\quad=\frac{K\omega_n 2^{n-2}}{n}\left(C_1\int_{\mathbb{R}^n}A(|\nabla u_k|)dy\right)^\mu. \int_{l_k^n}^\infty B\circ A^{-1}\left({{1}\over t}\right)dt 
<\infty
\end{eqnarray*}

\n It follows that $|\nabla u_k|\in L^B(\mathbb{R}^n)$. Therefore, (3.6) is valid for $u_k$ and we have
\begin{equation}\label{e3.10}
\int_{\mathbb{R}^n} B(|\nabla u_k|)dx
\leq \frac{\gamma' C_9}{\delta^{\kappa}(1-2\gamma' C_9\delta^{1+a_0})}\int_{\mathbb{R}^n}B(|F_k|)dx
\quad\forall k
\end{equation}

\n Since $F_k\rightarrow F$ strongly in $L^B(\mathbb{R}^n)$, we deduce from (3.10) that we have for some positive 
integer $k_0$
\begin{equation}\label{e3.11}
\int_{\mathbb{R}^n} B(|\nabla u_k|)dx
\leq \frac{2\gamma' C_9}{\delta^{\kappa}(1-2\gamma' C_9\delta^{1+a_0})}\int_{\mathbb{R}^n}B(|F|)dx
\quad\forall k\geq k_0
\end{equation}

\n Therefore, $u_k$ is uniformly bounded in $W^{1,B}(\mathbb{R}^n)$ and consequently has a weakly convergent
subsequence to some function $v$ in $W^{1,B}(\mathbb{R}^n)$. Using Lemma 2.7 for the function $B$, we get
\begin{equation}\label{e3.12}
\int_{\mathbb{R}^n} B(|\nabla u_k|)dx\geq \int_{\mathbb{R}^n} B(|\nabla v|)dx+\int_{\mathbb{R}^n}\left<\frac{b(|\nabla v|)}{|\nabla v|}\nabla v,\nabla u_k-\nabla v\right>dx
\end{equation}

\n Passing to the limit in (3.11) and taking into account (3.12) and the convergence for the weak topology, we get
\begin{equation}\label{e3.13}
\int_{\mathbb{R}^n} B(|\nabla v|)dx
\leq \frac{2\gamma' C_9}{\delta^{\kappa}(1-2\gamma' C_9\delta^{1+a_0})}\int_{\mathbb{R}^n}B(|F|)dx
\end{equation}

\n The proof of the theorem will be complete if we prove that $v$ is a solution of problem $(P)$.
To this end, we observe that it is enough to show that $\Theta(\nabla u_k)$ has a 
weakly convergent subsequence in $L^{\widetilde{A}}(\mathbb{R}^n)$ to $\Theta(\nabla v)$.
The proof is well known for monotone and continuous operators. We give it here for the sake of
completeness.

\n Using $w=u_k-v$ as a test function for problem $(P_k)$, we get 
\begin{eqnarray}\label{e3.14}
\int_{\mathbb{R}^n} \Theta(\nabla u_k).(\nabla u_k-\nabla v)dx
=\int_{\mathbb{R}^n}\Theta(F_k).(\nabla u_k-\nabla v) dx
\end{eqnarray}

\n Now, given that $\nabla u_k$ is uniformly bounded in $L^{A}(\mathbb{R}^n)$,
$\Theta(\nabla u_k)$ is uniformly bounded in
$L^{\widetilde{A}}(\mathbb{R}^n)$, and therefore has a weakly convergent
subsequence to some vector function $V$ in $L^{\widetilde{A}}(\mathbb{R}^n)$. 
Moreover, $u_k\rightharpoonup v$ weakly in $W^{1,B}(\mathbb{R}^n)$ 
and $F_k\rightarrow F$ strongly in $L^{B}(\mathbb{R}^n)$. So, in particular, 
$u_k\rightharpoonup v$ weakly in $W^{1,A}(\mathbb{R}^n)$
and $F_k\rightarrow F$ strongly in $L^{A}(\mathbb{R}^n)$.
Therefore, we obtain from (3.14)
\begin{eqnarray}\label{e3.15}
&&\limsup_{k\rightarrow\infty}\int_{\mathbb{R}^n}|\nabla u_k|a(|\nabla u_k|)dx=\int_{\mathbb{R}^n} V.\nabla vdx=
\end{eqnarray}

\n At this point, we use (1.7) for $X=\nabla u_k$ and an arbitrary vector function $Y$ in $L^{A}(\mathbb{R}^n)$
\begin{eqnarray*}
\int_{\mathbb{R}^n} \big(\Theta(\nabla u_k)-\Theta(Y)\big). (\nabla u_k-Y) dx\geq 0.
\end{eqnarray*}
Letting  $k\rightarrow\infty$ and taking into account (3.15), one can check that
\begin{eqnarray*}
\int_{\mathbb{R}^n}\big(V-\Theta(Y)\big). (\nabla v-Y) dx\geq0.
\end{eqnarray*}
Now choosing $Y=\nabla v-\lambda\vartheta$, where $\lambda$ is
a positive number and $\vartheta$ is an arbitrary vector
function that belongs to $\mathcal{D}(\mathbb{R}^n)$, we obtain
\begin{eqnarray*}
\int_{\mathbb{R}^n} \big(V-\Theta(\nabla v-\lambda\vartheta)\big).
\vartheta dx\geq 0.
\end{eqnarray*}
Letting $\lambda\rightarrow 0$, we get
$\displaystyle{\int_{\mathbb{R}^n} \big(V-\Theta(\nabla v)\big). \vartheta
dx\geq 0.}$ Since $\vartheta$ is an arbitrary vector function in 
$\mathcal{D}(\mathbb{R}^n)$, it follows that $V=\Theta(\nabla v)$ in $\mathbb{R}^n$. 
This achieves the proof of Theorem 1.1.
\qed

\end{document}